\documentclass[preprint,reqno,12pt]{elsarticle}

\usepackage{amssymb}
\usepackage{amsmath}
\usepackage{amsthm}
\usepackage{algpseudocode}
\usepackage[linesnumbered,ruled,vlined]{algorithm2e}
\usepackage[normalem]{ulem}
\usepackage{xcolor}

\newcommand{\smallfrac}[2]{\textstyle\frac{#1}{#2}\displaystyle}

\newcommand{\NN}{\mathbb{N}}

\newcommand{\finwords}{\{0,1\}^*}
\newcommand{\infwords}{\{0,1\}^{\omega}}
\newcommand{\words}[1]{\{0,1\}^{#1}}

\newcommand{\Sn}{\mathfrak{S}_n}
\newcommand{\Q}{\mathcal{Q}}
\newcommand{\A}{\mathcal{A}}
\renewcommand{\P}{\mathcal{P}_S(n)}
\DeclareMathOperator{\st}{st}
\DeclareMathOperator{\Des}{Des}
\DeclareMathOperator{\Pk}{Pk}
\DeclareMathOperator{\upgr}{\overline{gr}}
\DeclareMathOperator{\logr}{\underline{gr}}
\DeclareMathOperator{\gr}{gr}
\DeclareMathOperator{\Alt}{Alt}

\newtheorem{thm}{Theorem}[section]
\newtheorem{prop}[thm]{Proposition}
\newtheorem{lem}[thm]{Lemma}

\newtheorem{prob}[thm]{Problem}

\newtheorem{rmk}[thm]{Remark}

\numberwithin{equation}{section}

\journal{European Journal of Combinatorics}

\begin{document}

\begin{frontmatter}



\title{Growth Rates of Permutations with Given Descent or Peak Set}


\author{Mohamed Omar} %% Author name

\affiliation{organization={Department of Mathematics \& Statistics, York University},%Department and Organization
            addressline={4700 Keele St}, 
            city={Toronto},
            postcode={M3J 1P3}, 
            state={Ontario},
            country={Canada}}

\author{Justin M. Troyka} %% Author name

\affiliation{organization={Department of Mathematics, California State University, Los Angeles},%Department and Organization
            addressline={5151 State University Dr}, 
            city={Los Angeles},
            postcode={90032}, 
            state={California},
            country={United States}}

%% Abstract
\begin{abstract}
Given a set $I \subseteq \mathbb{N}$, consider the sequences $\{d_n(I)\},\{p_n(I)\}$ where for any $n$, $d_n(I)$ and $p_n(I)$ respectively count the number of permutations in the symmetric group $\Sn$ whose descent set (respectively peak set) is $I \cap [n-1]$. We investigate the growth rates $\gr d_n(I) = \lim_{n \to \infty} \left(d_n(I)/n!\right)^{1/n}$ and $\gr p_n(I) = \lim_{n \to \infty} \left(p_n(I)/n!\right)^{1/n}$ over all $I \subseteq \mathbb{N}$. Our main contributions are two-fold. Firstly, we prove that the numbers $\gr d_n(I)$ over all $I \subseteq \mathbb{N}$ are exactly the interval $\left[0,2/\pi\right]$. To do so, we construct an algorithm that explicitly builds $I$ for any desired limit $L$ in the interval. Secondly, we prove that the numbers $\gr p_n(I)$ for periodic sets $I \subseteq \mathbb{N}$ form a dense set in $\left[0,1/\sqrt[3]{3}\right]$. We do this by explicitly finding, for any prescribed $L$ in the interval, a set $I$ whose corresponding growth rate is arbitrarily close to $L$.
\end{abstract}



%% Keywords

\end{frontmatter}



\section{Introduction}

For a set $I \subseteq \mathbb{N}$, let $d_n(I)$ be the number of permutations in the symmetric group $\mathfrak{S}_n$ that have descents precisely at $I \cap [n-1]$; in other words, $d_n(I) := \#\{\pi \in \Sn \colon \pi_i>\pi_{i+1} \mbox{ if and only if } i \in I \cap [n-1] \}$. Understanding the statistic $d_n(I)$ has long been of interest to mathematicians. A classical result of this kind is enumerating \emph{alternating permutations} in $\Sn$: these are permutations $\pi_1\pi_2 \cdots \pi_n$ for which either
\[
\pi_1>\pi_2<\pi_3>\pi_4< \cdots \ \text{ or } \ \pi_1<\pi_2>\pi_3<\pi_4> \cdots
\]
and are therefore counted by $d_n(I)$ where respectively, $I=\{1,3,5,\ldots\}$ and $I=\{2,4,6,\ldots\}$. Regardless of which of these two sets $I$ we select, $d_n(I)$ is the same, observed through the bijection that sends $\pi \mapsto \pi'$ where $\pi_i'=(n+1)-\pi_i$ for each $i$. It is classically known that for these two sets $I$, $\{d_n(I)\}$ is the well-known combinatorial sequence $\{E_n\}$ called the Euler numbers, also known as the zigzag numbers or ``tangent and secant" numbers. It is well-known that these numbers have exponential generating function given by $\tan(x)+\sec(x)$. 

Outside of alternating permutations, the numbers $d_n(I)$ for $I \subseteq \mathbb{N}$ have been investigated in many contexts. MacMahon \cite{Macmahon} studied them in depth, deriving the following explicit formula:
\begin{thm} \label{thm:macmahon}
Let $n \ge 1$ be an integer. For a set $J  = \{j_1 < j_2 < \dots < j_r\} \subseteq [n-1]$, define $\operatorname{co}(J) = (j_1, j_2 - j_1, j_3 - j_2, \dots, j_r - j_{r-1}, n - j_r)$ to be the integer composition corresponding to $J$. Then, for each $I \subseteq [n-1]$,
\[ d_n(I) = \sum_{J \subseteq I} (-1)^{|I|-|J|} \binom{n}{\operatorname{co}(J)}, \]
where $\binom{n}{\mu_1, \dots, \mu_k} = \frac{n!}{\mu_1 ! \dots \mu_k !}$ denotes the multinomial coefficient.
\end{thm}
From this it follows that $d_n(I)$ is polynomial in $n$ if $I$ is finite. These polynomials were recently studied by Diaz-Lopez, Harris, Insko, Omar and Sagan in \cite{DHIOS} where they found combinatorial interpretations for the coefficients and explored roots to address questions of unimodality and log-concavity. Expanded work on this includes the work of Bencs \cite{Bencs} and references therein. When $I$ is infinite, many of the results in the literature establish asymptotic formulae in $n$ for $d_n(I)$. Elizalde and Troyka \cite[Sec.\ 5]{ET} show that $d_n(I)$ is greater than $(n/2) \lfloor n/2\rfloor !$ for ``almost all'' sets $I$. Bender, Helton and Richmond \cite{BHR} prove an asymptotic formula of the form $d_n(I) \sim cL^n n!$ when $I$ is ``nearly periodic''. Other prior work on this  includes the work of Brown, Fink and Willbrand \cite{BFW}, Luck \cite{Luck} and references therein. The latter, from the point of view of statistical physics, offers a significant alternative perspective on our work --- see Problem \ref{prob:distribution} and the remarks that follow it.

It is evident then that a focus in the literature is the asymptotics of $d_n(I)$ as $n$ grows, taking $I$ to range over all subsets of $\mathbb{N}$. An intriguing question then is what the range of the asymptotic behavior of $d_n(I)$ is. On one extreme end, Niven \cite{Niven} and de Bruijn \cite{deBruijn} independently proved that the sets $I$ that achieve the largest value of $d_n(I)$ are precisely those that give rise to alternating permutations. In fact they proved the following more general result:

\begin{thm} \label{thm:alternation-set} \cite[Thm.\ 5]{Niven} \cite[Sec.\ 3]{deBruijn} \cite[Prop.\ 1.6.4]{Stanley}
Fix $n \in \mathbb{N}$. For each $S \subseteq [n-1]$, define $\Alt(S) \subseteq [n-2]$ to be the \emph{alternation set} of $S$, i.e.\
\[ \Alt(S) = \{i \in [n-2] \colon \text{exactly one of $i$ and $i+1$ is in $S$}\}. \]
For all $S, T \subseteq [n-1]$, if $\Alt(S) \subsetneqq \Alt(T)$, then $d_n(S) < d_n(T)$. In particular, the two subsets $I \subseteq \mathbb{N}$ that maximize $d_n(I)$ are
\[
I=\{1,3,5,\ldots\} \ \ \text{ and } \ \ I=\{2,4,6,\ldots\}.
\]
\end{thm}

For both such sets $I$ which maximize $d_n(I)$, a quick asymptotic analysis of the  exponential generating series $\tan(x)+\sec(x)$ for $\{d_n(I)\} = \{E_n\}$ shows that the proportion of permutations in $\Sn$ with descent set $I \cap [n-1]$ is 
\[
\frac{d_n(I)}{n!} = \frac{4}{\pi} \left( \frac{2}{\pi} \right)^n + O\left( \left(\frac{2}{3\pi}\right)^n\right).
\]

Consequently, the proportion is asymptotically $(4/\pi) \cdot (2/\pi)^n$. Since this statistic exponentially decays as a function of $n$, it is often easier to reduce it to a single number, so we retain the statistic $\lim_{n \to \infty} \left(d_n(I)/n!\right)^{1/n} = 2/\pi$.  For any $I$, the limit $\lim_{n \to \infty} \left(d_n(I)/n!\right)^{1/n}$ (when it exists) is called the \emph{growth rate} of $\{d_n(I)\}$ and is denoted by $\gr d_n(I)$. Growth rates are not only convenient but are well-studied for sequences arising in permutation enumeration problems. A ubiquitous example of this is the study of Stanley--Wilf limits of permutations, which classify growth rates of permutation classes that avoid a given permutation pattern (see the survey \cite{vatter2015permutation} for more). We see then that for any $I \subseteq \mathbb{N}$, $\gr d_n(I) \leq 2/\pi$. In the other extreme, again by the work of Diaz-Lopez et al.\@ \cite{DHIOS}, if $I$ is finite, then $d_n(I)$ is polynomial in $n$ and so $\gr d_n(I) = 0$. This gives us a range of possible values for the growth rate $\gr d_n(I)$ as $I$ varies:
\[
0 \leq \gr d_n(I) \leq 2/\pi.
\]
A natural question arises: for which $L \in \left[0,2/\pi\right]$ is $L$ the growth rate  $\gr d_n(I)$ for some $I \subseteq \mathbb{N}$? Our first main theorem shows this happens for \emph{every} $L$.

\begin{thm} \label{thm:properexists}
For every $L \in [0, 2/\pi]$, there exists $I \subseteq \mathbb{N}$ such that the growth rate $\gr d_n(I)$ exists and equals $L$.
\end{thm}

For many sets $I \subseteq \mathbb{N}$, the sequence $(d_n(I)/n!)^{1/n}$ does not converge. But even so, one can still investigate the limiting behavior of this sequence by looking at the \emph{lower growth rate} $\logr d_n(I) := \liminf_{n \to \infty} (d_n(I)/n!)^{1/n}$ and the \emph{upper growth rate} $\upgr d_n(I) := \limsup_{n \to \infty} (d_n(I)/n!)^{1/n}$, both of which always exist. Taking this into consideration, we obtain a strengthening of Theorem~\ref{thm:properexists} (see Theorem~\ref{thm:infsup}) that shows that if $0 \leq L \leq L' \leq 2/\pi$ then we can find a set $I$ such that $\logr d_n(I)=L$ and $\upgr d_n(I)=L'$.

Many enumerative results pertaining to descent sets of permutations have analogues for peak sets. For $I \subseteq \mathbb{N}$, let $p_n(I)$ be the number of permutations in $\Sn$ that have peaks precisely at $I \cap [n-1]$, that is $p_n(I) := \# \{\pi \in \Sn \colon \pi_{i-1} < \pi_i > \pi_{i+1} \text{ if and only if } i \in I \cap [n-1]\}$. Note by our convention, a permutation never has a peak at its first and last position (however if one wants to add such a stipulation, many of our results can be adapted to accommodate accordingly). Similar to the case for descents, when $I$ is finite it was proven in \cite{DHIO} that $p_n(I)$ is polynomial in $n$ multiplied by a power of $2$ that is linear in $n$ and $|I|$. More on these polynomials can be found in \cite{DHIO} and \cite{KantarciOguz}. Again, analogous to the case for descents, Kasraoui proved that $p_n(I)$ is maximized for the following sets:

\begin{thm}\cite[Thm.\ 6.1]{BFT}
The set $I$ that maximizes $p_n(I)$ is

\[
I=
\begin{cases}
\{3,6,9,\ldots\} \cap [n-1] \text{ and } \{4,7,10,\ldots,\} &\mbox{ if } n \equiv 0 \ (\text{mod } 3), \\
\{3,6,9,\ldots,3s,3s+2,3s+5,\ldots\} \cap [n-1] & \mbox{ if } n \equiv 1 \ (\text{mod } 3), \\
\{3,6,9,\ldots\} \cap [n-1] & \mbox{ if } n \equiv 2 \ (\text{mod } 3). \\
\end{cases}
\]
\end{thm}
Here, when $n \equiv 1 \ (\text{mod } 3)$, $s$ is any integer with $1 \leq s \leq \left \lfloor \frac{n}{3} \right \rfloor$. Kasraoui further enumerates $p_n(I)$ for each of these sets $I$ and in all cases one can extrapolate that the growth rate $\lim_{n \to \infty} \left(p_n(I)/n!\right)^{1/n}$, which we denote by $\gr p_n(I)$, is exactly $1/\sqrt[3]{3}$. A similar question to the one we asked for descents now arises: for which $L \in \left[0,1/\sqrt[3]{3}\right]$ is $L$ the growth rate $\gr p_n(I)$ of some $I \subseteq \mathbb{N}$? Our main contribution here is that, even looking only at periodic sets $I$, growth rates are dense in the interval:

\begin{thm}\label{thm:peakgrowthrate}
For any $L \in  \left[0,1/\sqrt[3]{3}\right]$ and $\epsilon>0$, there exists a periodic set $I \subseteq \mathbb{N}$ such that the growth rate $\gr p_n(I)$ satisfies
\[
\left|\gr p_n(I) - L\right|<\epsilon.
\]
\end{thm}

Our paper is organized as follows. In order to make our results transparent it will be much more convenient to work in the domain of combinatorics on words. Section~\ref{sec:def} introduces all necessary definitions and notation that translate concepts on descents and peaks into this framework. The main contributions of the article are in Sections~\ref{sec:descents} and \ref{sec:peaks}, with Section~\ref{sec:descents} dedicated to the proof of Theorem~\ref{thm:properexists} and Section~\ref{sec:peaks} dedicated to the proof of Theorem~\ref{thm:peakgrowthrate}. We end in Section~\ref{sec:conclusion} with future directions and open problems.

\section{Definitions and notation}\label{sec:def}

Throughout our exposition it will be easier to discuss descents and peaks at positions $I \subseteq \mathbb{N}$ by referring to binary strings that have a $1$ in positions where we have a descent (or respectively peak) and $0$ otherwise. This motivates us to introduce the language of words relevant to our context. Throughout this article, we will work separately with finite and infinite words in the alphabet $\{0,1\}$.

Let $\words{k}$ denote the set of words of length $k$, i.e.\@ $\words{k} = \{x_1 \cdots x_k \mid x_1, \ldots, x_k \in \{0,1\}\}$. Let $\finwords = \bigcup_{k\ge0} \words{k}$ denote the set of finite $\{0,1\}$ words. In the case we find the finiteness condition too limiting, we use $\infwords$ to denote the set of $\{0,1\}$ words indexed by the positive integers, i.e.\@ $$\infwords = \{x_1 x_2 \cdots \mid x_1, x_2, \ldots \in \{0,1\}\}.$$ Our paradigm will often have us constructing words with repetition in them. In that light, if $x$ is a word, then $x^r$ denotes the $r$-fold concatenation of $x$ with itself: that is, $x^r = x \cdots x$ with $r$ copies of $x$. There are situations in which we may concatenate a word to itself a countably infinite number of times, so for a word $x$, we let $x^\omega$ denote the infinite concatenation of $x$ with itself, i.e.\@ $x^{\omega} = xx\cdots$. For odd integers $r$, it will be convenient to slightly abuse notation here: we write  $(01)^{r/2}$ to denote $(01)^{(r-1)/2}0$, and similarly we write $(10)^{r/2}$ for $(10)^{(r-1)/2}1$. The intent here is that both $(01)^{r/2}$ and $(10)^{r/2}$ remain strings of length $r$. If $x$ is a $\{0,1\}$-word whose length is at least $k$, let $x|_{k} \in \words{k}$ denote the word formed by the first $k$ letters of $x$; that is, if $x = x_1 \cdots x_j \in \words{j}$ for $j \ge k$ or $x = x_1 x_2 \cdots \in \infwords$, then $x|_k = x_1 \cdots x_k$. As an example to clarify notation here, if $x=(10)^8$ then $x|_5=10101=(10)^{5/2}$ using our convention.

We now contextualize the combinatorics on words framework by bridging it into our understanding of peaks and descents. Given $\pi \in \Sn$, the \emph{descent word} of $\pi$ is the word $\Des(\pi) = w_1 \cdots w_{n-1} \in \words{n}$ where $w_i = 1$ if $\pi_i > \pi_{i+1}$ (i.e. $i$ is a \emph{descent} of $\pi$) and $w_i = 0$ otherwise. In a similar light, the \emph{peak word} of $\pi$ is the word $\Pk(\pi) = w_1 \cdots w_{n-1} \in \words{n-1}$ where $w_i = 1$ if $\pi_{i-1} < \pi_i > \pi_{i+1}$ ($i$ is a \emph{peak} of $\pi$) and $w_i = 0$ otherwise. Note that $w_1 = 0$ for a peak word since we cannot have a peak at the first entry. For instance, if $\pi = 31248576 \in \mathfrak{S}_8$ then $\Des(\pi)=1000101$ whereas $\Pk(\pi)=0000101$. Given a $\{0,1\}$-word $w$ whose length is at least $n-1$ we can count the number of permutations whose descent or peak word matches the  appropriately sized prefix of $w$. Define $d_n(w)$ to be the number of permutations in $\mathfrak{S}_n$ with descent word $w|_{n-1}$; that is, $$d_n(w) = \# \{\pi \in \mathfrak{S}_n \colon \Des(\pi) = w|_{n-1}\}.$$ Similarly, if $w$ has length at least $n-1$, define $p_n(w)$ to be the number of permutations in $\mathfrak{S}_n$ with peak word $w|_{n-1}$; that is, $$p_n(w) = \# \{\pi \in \mathfrak{S}_n \colon \Pk(\pi) = w|_{n-1}\}.$$ If $\Des(\pi) = w|_{n-1}$, we will sometimes say that $\pi \in \mathfrak{S}_n$ has descent word $w$, even when the length of $w$ is greater than $n-1$; we may use similar language for the peak word of $\pi$. 

Notice that if $I \subseteq \mathbb{N}$, then $d_n(I)=d_n(w)$ and $p_n(I)=p_n(w)$ where $w$ is the word with $w_i=1$ if and only if $i \in I$. In that light then, our central questions in this article can be cast in terms of the sequences $\{d_n(w)\}$ and $\{p_n(w)\}$ for fixed $w \in \infwords$. This is precisely the approach used in Sections~\ref{sec:descents} and \ref{sec:peaks}.

\section{Descents}\label{sec:descents}

This section is dedicated to the proof of Theorem~\ref{thm:properexists}, which in the language of words states that every number in $[0, 2/\pi]$ is the growth rate $\gr d_n(w)$ of some $w \in \infwords$. Our approach necessitates proving identities and inequalities about the numbers $d_n(w)$, including the fact that two infinite words $w, w'$ with the same tail have the same growth rate; i.e. $\gr d_n(w) = \gr d_n(w')$. These results may be of interest in their own right, but their main purpose is to serve the main theorem, whose proof relies on an algorithm we develop that constructs an infinite word with repeating blocks of $0$'s and $10$'s in the correct proportions to achieve a given growth rate $L \in [0, 2/\pi]$.

In some of our proofs, we make use of a convenient construct for permutations: given $\alpha \in \NN^k$ whose entries are distinct, the \emph{standardization} of $\alpha$, denoted $\st(\alpha)$, is the permutation in $\mathfrak{S}_k$ with the same relative order as $\alpha$ --- more precisely, $\st(\alpha)$ is the permutation obtained from $\alpha$ by replacing the $i$th lowest number in $\alpha$ with $i$ for all $i$. Note that the list of integers $\alpha$ and $\st(\alpha)$ appear in the same relative order.

To start, we prove Lemma \ref{lem:binomialsum} which generalizes \cite[Prop.\@ 2.1]{DHIOS}, a proposition relating the enumeration of permutations with a given descent set to those with the same descent set but one element removed. Then, we iterate Lemma \ref{lem:binomialsum} to obtain Lemma \ref{lem:descent-formula}, a convenient recursive formula for the numbers $d_n(w)$.

\begin{lem} \label{lem:binomialsum}
Let $a,b \geq 1$ be integers. If $u \in \words{a-1}$ and $v \in \words{b-1}$, then
\[ \binom{a+b}{a} d_a(u) d_b(v) = d_{a+b}(u0v) + d_{a+b}(u1v). \]
\end{lem}

\begin{proof}
The right side precisely counts the set
\[ X = \{\pi \in \mathfrak{S}_{a+b} \colon  \text{$\Des(\pi) = u0v$ or $\Des(\pi) = u1v$}\}. \]
Now, $X$ is the set consisting of permutations $\pi = \pi_1 \cdots \pi_{a+b} \in \mathfrak{S}_{a+b}$ such that $\st(\pi_1 \cdots \pi_a)$ has descent word $u$ and $\st(\pi_{a+1} \cdots \pi_b)$ has descent word $v$. There are $\binom{a+b}{a}$ ways to choose which values appear in $\pi_1 \cdots \pi_a$, after which there are $d_a(u)$ ways to choose $\pi_1 \cdots \pi_a$ and $d_b(v)$ ways to choose $\pi_{a+1} \cdots \pi_{a+b}$. Thus,
\begin{equation*} %\label{eq:secondcount}
|X| = \binom{a+b}{a}d_a(u)d_b(v). \qedhere
\end{equation*}
\end{proof}
As an example to illustrate Lemma~\ref{lem:binomialsum}, let $a=6$, $b=3$, $u=10010 \in \{0,1\}^5$, and $v=00 \in \{0,1\}^2$. In this case, Lemma~\ref{lem:binomialsum} yields
\[ \binom{9}{6} d_6(10010) d_3(00) = d_9(10010000) + d_9(10010100). \]
Using the fact that $d_3(00) = 1$ (only one permutation in $\mathfrak{S}_3$ has no descents), we can rearrange the equation to obtain
\[
d_9(10010100)=\binom{9}{6} d_6(10010) - d_9(10010000).
\]
This illustrates how we can use Lemma~\ref{lem:binomialsum} to express $d_n(w)$ in terms of words with one fewer $1$ in them. Iterating this process yields the following lemma:

\begin{lem} \label{lem:descent-formula}
Let $w \in \infwords$ and $n \ge 1$. If the $1$'s in $w|_{n-1}$ are in positions $i_1 <  \cdots < i_k$, then
\[ d_n(w) = \sum_{r=0}^k (-1)^{k-r} \binom{n}{i_r} d_{i_r}(w), \]
where by convention we set $i_0 = 0$ and $d_0(w) = 1$.
\end{lem}

\begin{proof}
Let $r$ be an integer, $1 \le r \le k$. We apply Lemma \ref{lem:binomialsum} with $u = w_1 \cdots w_{i_r-1}$ and $v = 0^{n-i_r-1}$: using the fact that $d_{i_r}(w_1 \cdots w_{i_r-1}) = d_{i_r}(w)$ and $d_{n-i_r}(0^{n-i_r-1}) = 1$ yields
\[  d_n(w_1 \cdots w_{i_r} 0^\omega)=\binom{n}{i_r} d_{i_r}(w) - d_n(w_1 \cdots w_{i_{r-1}} 0^\omega). \]
Now iteratively apply Lemma~\ref{lem:binomialsum} while updating $u$ and $v$ appropriately at each step to get
\begin{align*}
d_n(w_1 \cdots w_{i_r} 0^\omega)&=\binom{n}{i_r} d_{i_r}(w) - d_n(w_1 \cdots w_{i_{r-1}} 0^\omega) \\
&=\binom{n}{i_r} d_{i_r}(w) \\& - \left( \binom{n}{i_{r-1}} d_{i_{r-1}}(w)\underbrace{d_{n-i_{r-1}}(0^{n-i_{r-1}-1})}_{=1} - d_n(w_1 \cdots w_{i_{r-2}}0^\omega) \right) \\
&=\binom{n}{i_r} d_{i_r}(w) - \binom{n}{i_{r-1}} d_{i_{r-1}}(w) + d_n(w_1 \cdots w_{i_{r-2}}0^\omega) \\
& \ \vdots \\
&=\sum_{r=0}^k (-1)^{k-r} \binom{n}{i_r} d_{i_r}(w).
\end{align*}
\end{proof}
Iteratively applying Lemma \ref{lem:descent-formula} gives rise to a very similar formula to the classically known explicit formula we stated as Theorem \ref{thm:macmahon}.
\begin{prop}
Let $w \in \{0,1\}^\omega$ and $n \geq 1$. If the $1$'s in $w|_{n-1}$ are in position $i_1<\cdots<i_k$ then
\[
d_n(w)=\sum_{J \subseteq \{i_1,\ldots,i_k\}} (-1)^{k-|J|} \binom{n}{\delta(J)}
\]
where here, if $J=\{j_1,j_2,\ldots,j_{\ell}\}$ then $\binom{n}{\delta(J)} = \frac{n!}{j_1!(j_2-j_1)! \cdots (j_{\ell}-j_{\ell-1})!(n-j_{\ell})!}$.
\end{prop}

\subsection{Inequalities, the shift operator, and words with the same tail}

In this subsection we prove a number of inequalities which will help us understand how the statistic $d_n(w)$ changes when we shift $w$, culminating in our result that words with the same tail (which we call \emph{equicaudal}) have the same growth rate (Proposition~\ref{prop:equicaudal}). The inequalities in this subsection will also be used in the proof of our main theorem on the existence of a word that achieves any specified growth rate (Theorem \ref{thm:properexists}).

\begin{lem} \label{lem:another}
If $w \in \infwords$ and $n \ge 1$, then
\[ d_{n-1}(w) \le d_n(w) \le (n-1) d_{n-1}(w). \]
\end{lem}

\begin{proof}
Given a permutation $\pi \in \mathfrak{S}_{n-1}$ with descent word $w$, we can form a permutation in $\Sn$ by choosing $i \in [n]$, increasing by $1$ every value in $\pi$ which is at least $i$, and appending $i$ to the end of $\pi$. There are $n$ options for $i$, and the number of those options which result in a permutation with descent word $w$ is between $1$ and $n-1$.
\end{proof}

\begin{lem} \label{lem:binomialinequality}
Let $n,k \ge 1$ and $w \in \words{n-1}$. Then
\begin{itemize}
\item[(a)] $d_{n+k}(w 0^k) \le \binom{n+k}{k} d_n(w)$.
\item[(b)] If $w$ ends with $0$, then $d_{n+k}(w (10)^{k/2}) \ge \frac12 \binom{n+k}{k} E_k d_n(w)$ (recall the convention that $(10)^{k/2}$ denotes $w(10)^{(k-1)/2} 1$ if $k$ is odd). Here, $E_k$ is the $k$th Euler number.
\end{itemize}
\end{lem}

\begin{proof}
(a) Apply Lemma \ref{lem:binomialsum} with $u = w$ and $v = 0^{k-1}$:
\[ \binom{n+k}{k} d_n(w) = d_{n+k}(w0^k) + d_{n+k}(w10^{k-1}) \ge d_{n+k}(w0^k). \]

(b) Write $w = z0$. Apply Lemma \ref{lem:binomialsum} with $u = z0$ and $v = (01)^{(k-1)/2}$: using the fact that $d_k( (01)^{(k-1)/2}) = E_k$ yields
\[ \binom{n+k}{k} d_n(w) E_k = d_{n+k}(z01(01)^{(k-1)/2}) + d_{n+k} (z00(01)^{(k-1)/2}). \]
But $d_{n+k} (z00(01)^{(k-1)/2}) \le d_{n+k}(z01(01)^{(k-1)/2})$, because the alternation set (the set of positions of occurrences of $10$ or $01$) in $z00(01)^{(k-1)/2}$ is contained in the alternation set of $z01(01)^{(k-1)/2}$ (see Theorem \ref{thm:alternation-set}). Thus
\[ \binom{n+k}{k} d_n(w) E_k \le 2d_{n+k}(z01(01)^{(k-1)/2}) = 2d_{n+k}(w (10)^{k/2}), \]
and dividing by $2$ gives the desired result.
\end{proof}

For $m \in \NN$ define the \emph{$m$-shift} of $w$ by
\[ \Sigma^m(w) = \Sigma^m(w_1 w_2 \cdots) = w_{m+1} w_{m+2} \cdots. \]
When $m=1$ we simply write the $1$-shift as $\Sigma(w)$. We say that two infinite words $w, w' \in \infwords$ are \emph{equicaudal} if there are finite words $v,v' \in \finwords$ (possibly of different lengths) and an infinite word $z \in \infwords$ such that $w = vz$ and $w' = v'z$; that is, $w$ and $w'$ have the same tail or the same coda. Note that two words $w, w' \in \infwords$ are equicaudal if and only if $\Sigma^m(w) = \Sigma^{m'}(w')$ for some $m, m' \in \NN$. Essential to our main theorem is the proof that equicaudal words have the same growth rates.

\begin{prop} \label{prop:equicaudal}
Let $w, w' \in \infwords$, If $\Sigma^m(w) = \Sigma^{m'}(w')$ for some $m, m' \in \NN$, then $\upgr d_n(w) = \upgr d_n(w')$ and $\logr d_n(w) = \logr d_n(w')$. Equivalently, equicaudal descent words have the same upper growth rate and the same lower growth rate.
\end{prop}

The proof of Proposition~\ref{prop:equicaudal} hinges on the following lemma, which can be used to prove that the upper and lower growth rates of a word are the same as the upper and lower growth rates of any tail of it.

\begin{lem} \label{lem:sigmalemma}
If $w \in \infwords$ and $n \ge 1$, then
\[ \smallfrac{1}{n} d_n(w) \le d_n(\Sigma(w)) \le n d_n(w). \]
\end{lem}

\begin{proof}
Let $A$ be the set of $\pi \in \Sn$ with descent word $w$, and let $B$ be the set of $\pi \in \mathfrak{S}_{n-1}$ with descent word $\Sigma(w)$. Define $\phi \colon A \to [n] \times B$ by
\[ \phi(\pi_1 \pi_2 \cdots \pi_n) = (\pi_1, \operatorname{st}(\pi_2 \cdots \pi_n)), \]
where $\pi_1 \pi_2 \cdots \pi_n$ denotes a permutation in one-line notation. Since $\phi$ is an injective function, and since $|A| = d_n(w)$ and $|B| = d_{n-1}(\Sigma(w))$, we obtain
\[ d_n(w) \le n d_{n-1}(\Sigma(w)). \]
Therefore,
\[ \smallfrac{1}{n} d_n(w) \le d_{n-1}(\Sigma(w)) \le d_n(\Sigma(w)) \]
(using Lemma \ref{lem:another}), proving the first inequality.

Now let $C$ be the set of $\pi \in \Sn$ with descent word $\Sigma(w)$, and let $D$ be the set of $\pi \in \mathfrak{S}_{n+1}$ with descent word $w$. If $w_1 = 0$ (so permutations in $D$ begin with an ascent), then define $\psi \colon C \to D$ by
\[ \psi(\pi_1 \cdots \pi_n) = \st(0 \pi_1 \cdots \pi_n); \]
similarly, if $w_1 = 1$ (so permutations in $D$ begin with a descent), then define $\psi \colon C \to D \times [n]$ by
\[ \psi(\pi_1 \cdots \pi_n) = \st((n+1) \pi_1 \cdots \pi_n). \]
Since $\psi$ is an injective function, and since $|C| = d_n(\Sigma(w))$ and $|D| = d_{n+1}(w)$, we obtain
\[ d_n(\Sigma(w)) \le d_{n+1}(w). \]
Therefore,
\[ d_n(\Sigma(w)) \le d_{n+1}(w) \le n d_n(w) \]
(using Lemma \ref{lem:another}), proving the second inequality.
\end{proof}

\begin{proof}[Proof of Proposition \ref{prop:equicaudal}] Let $w, w' \in \infwords$, and assume $\Sigma^m(w) = \Sigma^{m'}(w') = z$. By repeated application of Lemma \ref{lem:sigmalemma}, $n^{-m} d_n(w) \le d_n(z) \le n^m d_n(w)$. From the fact that $d_n(z) \le n^m d_n(w)$, we find that
\begin{align*}
\upgr d_n(z) &= \limsup_{n\to\infty} {\left( \frac{d_n(z)}{n!} \right)}^{1/n} \\
&\le \limsup_{n\to\infty} {\left( \frac{n^m d_n(w)}{n!} \right)}^{1/n} = \left(\lim_{n\to\infty} n^{m/n} \right) \upgr d_n(w) =  \upgr d_n(w);
\end{align*}
thus $\upgr d_n(z) \le \upgr d_n(w)$. From the fact that $d_n(z) \ge n^{-m} d_n(w)$, the same reasoning shows that $\upgr d_n(z) \ge d_n(w)$. Therefore, $\upgr d_n(z) = \upgr d_n(w)$.

Similarly, $\upgr d_n(z) = \upgr d_n(w')$. Therefore, $\upgr d_n(w) = \upgr d_n(w')$. Finally, the same reasoning with $\liminf_{n\to\infty}$ shows that $\logr d_n(w) = \logr d_n(w')$.
\end{proof}

\begin{rmk} \label{rmk:chirping}
An immediate corollary of Proposition \ref{prop:equicaudal} is that $\gr d_n(w) = 0$ if $w$ has only finitely many $1$'s (resp.\@ only finitely many $0$'s), since in this case $w$ and $0^\omega$ (resp.\@ $1^\omega$) are equicaudal. Interestingly, the converse is false: it is possible to have $\gr d_n(w) = 0$ even when $w$ has infinitely many $0$'s and infinitely many $1$'s. Indeed, Luck \cite[Sec.\@ 10]{Luck} analyzes ``chirping'' patterns of descents, in which the descents become more and more scarce as the binary word continues; he finds that, if $w$ is such a word, then $d_n(w)/n!$ is smaller than any exponential decay, and so $\gr d_n(w) = 0$. For instance, he finds that, if $w = 1(0^2) 1(0^4) 1 (0^6) 1 (0^8)\cdots$ is the infinite word with $1$'s in positions $r^2$ for $r\ge1$ and $0$'s in all other positions, then
\[ \ln \left(\frac{d_n(w)}{n!}\right) \sim -\frac{n}{2} (\ln(4n) - 3), \]
from which it follows that ${\left( \frac{d_n(w)}{n!} \right)\!}^{1/n} \le C n^{-p}$ for some constants $C,p>0$.
\end{rmk}

However, even though the growth rate can be $0$, the sequence $d_n(w)$ still grows faster than any polynomial if $w$ has infinitely many $0$'s and infinitely many $1$'s:

\begin{prop}
If $w$ has infinitely many $0$'s and infinitely many $1$'s, then, for every $p>0$, $d_n(w)$ is eventually greater than $n^p$.
\end{prop}

\begin{proof}
Let $p>0$. Since $w$ has infinitely many $0$'s and infinitely many $1$'s, there exists $i>p$ such that $w_i w_{i+1} = 10$. By Theorem \ref{thm:alternation-set}, we have
\[ d_n(w) \ge d_n(w|_i \,0^\omega), \]
since every position that has an alternation ($01$ or $10$) in $w|_i\, 0^\omega$ also has an alternation in $w$. By Lemma \ref{lem:descent-formula},
\[ d_n(w|_i \, 0^\omega) \sim \frac{n^i}{i!}, \]
which is eventually greater than $n^p$.
\end{proof}

\subsection{Existence of words with every possible growth rate}

The goal of this subsection is to prove the first main theorem, Theorem~\ref{thm:properexists}. We note that $\gr d_n(0^{\omega})=0$ and by our previous discussions $\gr d_n((01)^{\omega})=2/\pi$ so we may assume $L \in (0, 2/\pi)$. It will be convenient to also track $M=L^{-1}$.

To prove the theorem, we construct $w \in \infwords$ with $\gr d_n(w) = L$ using a recursive algorithm. The basic idea is as follows: while the algorithm is in ``state $0$'', we append $0$ repeatedly to decrease the growth rate; while the algorithm is in ``state $10$'', we  append $10$ repeatedly to increase the growth rate. To track the algorithm's state, we use a variable $q$ which is either $0$ or $10$ and we append $q$ to the right of the word at each step depending on constants $r_i$ that will serve as continuous indicators on whether to flip states. Our algorithm is as follows:

\begin{algorithm}[H]
 \KwData{$L \in \left(0,2/\pi\right), \ M=L^{-1}$}
 \KwResult{word $w \in \infwords$, sequence $\{r_n\}$}
 $q=0$, $n=1$, $r_1 = M$\;
 \While{}{
  \If{$r_n<1$}{
   $q \gets 10$\;}
   \If{$r_n>1$}{
   $q \gets 0$\;}
   \eIf{$q=0$}{
   $w_n \gets 0;$ \\
   $r_{n+1} \gets \frac{d_{n+1}(w|_n)}{(n+1)!}M^{n+1};$ \\
   $n \gets n+1$\;}{
   $w_n \gets 1, \ w_{n+1} \gets 0;$ \\
   $r_{n+1} \gets \frac{d_{n+1}(w|_n)}{(n+1)!}M^{n+1}, \ r_{n+2} \gets \frac{d_{n+2}(w|_{n+1})}{(n+2)!}M^{n+2};$ \\
   $n \gets n+2$\;}
   }
 \caption{Constructing Word with Growth Rate $L$}
\end{algorithm}
Algorithm 1 outputs an infinite word $w \in \infwords$ and a sequence of real numbers $\{r_n\} = \left \{\frac{d_n(w)}{n!} M^n \right \}$.

%\begin{itemize}
%\item[(1)] Check whether to update $q$:
%\begin{itemize}
%\item If $q = 0$ and $r_n < 1$, then set $q = 10$.
%\item If $q = 10$ and $r_n > 1$, then set $q = 0$.
%\item Else, leave $q$ unchanged.
%\end{itemize}
%\item[(2)] Next, append $0$ or $10$ to $w|_{n-1}$ depending on the state $q$:
%\begin{itemize}
%\item If $q = 0$, then set $w_n = 0$.
%\item If $q = 10$, then set $w_n w_{n+1} = 10$.
%\end{itemize}
%\end{itemize}

As an example, let $L = 1/2$, so $M = 2$. Algorithm 1 computes the following values:
\[ \begin{array}{l|l|l|l}
& q = 0 & & r_1 = 2 > 1 \\ \hline
n = 1 & q = 0 & w_1 = 0 & r_2 = \frac{d_2(0)}{2!} 2^2 = 2 > 1 \\ \hline
n = 2 & q = 0 & w_2 = 0 & r_3 = \frac{d_3(00)}{3!} 2^3 \approx 1.33 > 1 \\ \hline
n = 3 & q = 0 & w_3 = 0 & r_4 = \frac{d_4(000)}{4!} 2^4 \approx 0.67 < 1 \\ \hline
n = 4 & q = 10 & w_4 w_5 = 10 & r_5 = \frac{d_5(0001)}{5!} 2^5 \approx 1.07 \\
      &        &              & r_6 = \frac{d_6(00010)}{6!} 2^6 \approx 1.24 > 1 \\ \hline
n = 6 & q = 0 & w_6 = 0 & r_7 = \frac{d_7(000100)}{7!} 2^7 \approx 0.86 < 1 \\ \hline
n = 7 & q = 10& w_7 w_8 = 10 & r_8 = \frac{d_8(0001001)}{8!} 2^8 \approx 1.29 \\
     &       &              & r_9 = \frac{d_9(00010010)}{9!} 2^9 \approx 1.55 > 1 \\ \hline
n = 9 & q = 0 & w_9 = 0 & r_{10} = \frac{d_{10}(000100100)}{10!} 2^{10} \approx 1.09 > 1 \\ \hline
n = 10 & q = 0 & w_{10} = 0 & r_{11} = \frac{d_{11}(0001001000)}{11!} 2^{11} \approx 0.56 < 1 \\ \hline
n = 11 & q = 10& w_{11} w_{12} = 10 & r_{12} = \frac{d_{12}(00010010001)}{12!} 2^{12} \approx 0.89 \\
     &       &              & r_{13} = \frac{d_{13}(000100100010)}{13!} 2^{13} \approx 1.04 > 1 \\ \hline
\end{array} \]
Reading the outputs of $w_1, w_2, \dots$, we find that the word $w$ begins as $w = 000100100010 \dots$.

\begin{rmk}
We clarify the infinite nature of Algorithm 1. If we stop the algorithm after finitely many iterations, then the output is a finite word $w|_n = w_1 \dots w_n$ and a finite list of real numbers $r_1, \dots, r_n$ (for some $n$). But by Theorem \ref{prop:equicaudal}, we can replace $w_1 \dots w_n$ with a completely different word $w'_1 \dots w'_n$ leading to a different finite list of real numbers $r_1',\ldots,r_n'$, but still end up with the same growth rate. That is, no matter how long we run the algorithm before stopping it, we can erase all the output produced so far without affecting the growth rate. The output after a finite stopping point does not matter --- only the output's ``end behavior'' as $n \to \infty$. This is analogous to the fact that the first finitely many terms of an infinite series do not affect whether the series converges.

If we want to actually run the algorithm on a computer, then the given number $L$ must be a computable real number, of which there are only countably many. If $L$ is not computable, then it cannot be represented in a finite way for the computer to read as input. In order to prove Theorem \ref{thm:properexists} for all real numbers $L \in (0, 2/\pi)$, we must understand our algorithm not as a practical computer program, but as a recursive definition of a function $f \colon (0, 2/\pi) \to \infwords \times \mathbb{R}^\omega$ (where $\mathbb{R}^\omega$ denotes the set of sequences of real numbers). The recursion theorem from set theory tells us that this $f$ exists and is uniquely determined by its recursive definition. Thus, for all $L \in (0, 2/\pi)$, there exists $f(L) = (w, \{r_n\})$ --- even if $L$ is non-computable.
\end{rmk}

\begin{rmk} \label{rmk:suffices-to-show}
We have $\lim_{n\to\infty} {r_n}^{1/n} = M \gr d_n(w) = L^{-1} \gr d_n(w)$; so, to prove the theorem, it suffices to show that $\lim_{n\to\infty} {r_n}^{1/n} = 1$. We will do this by proving a few technical lemmas about $r_n$.
\end{rmk}

\begin{lem} \label{lem:r_n-inequality}
For all $n$,
\[ \frac{M}{n} r_{n-1} \le r_n \le M r_{n-1}. \]
\end{lem}

\begin{proof}
from Lemma \ref{lem:another}, we have $d_n(w) \ge d_{n-1}(w)$. So
\[ r_n = \frac{d_n(w)}{n!} M^n \ge \frac{d_{n-1}(w)}{n!} M^n = \frac{r_{n-1}}{n} M, \]
proving the first inequality. Also from Lemma \ref{lem:another}, we have $d_n(w) \le n d_{n-1}(w)$. So
\[ r_n = \frac{d_n(w)}{n!} M^n \le \frac{n d_{n-1}(w)}{n!} M^n = \frac{d_{n-1}(w)}{(n-1)!} M^n = r_{n-1} M, \]
proving the second inequality.
\end{proof}

\begin{lem} \label{lem:inequality-induction}
For $n,k\ge1$,
\[ \frac{M^k}{(n+k)^k} r_n \le r_{n+k} \le M^k r_n. \]
\end{lem}
\begin{proof}
Induction on $k$, using Lemma \ref{lem:r_n-inequality}.
\end{proof}

\begin{lem} \label{lem:0^k-inequality}
For $n,k \geq 1$:
\begin{itemize}
\item[(a)] If $w_n w_{n+1} \cdots w_{n+k-1} = 0^k$, and if $k$ is large enough that $M^k/k! \le 1$, then $r_{n+k} \le r_n$.
\item[(b)] If $w_n w_{n+1} \cdots w_{n+k-1} = (10)^{k/2}$ (recall our convention that $(10)^{k/2}$ denotes $(10)^{(k-1)/2}1$ if $k$ is odd), and if $k$ is large enough that $\frac{1}{2} \cdot \frac{E_k}{k!} \cdot M^k \ge 1$, then $r_{n+k} \ge r_n$.
\end{itemize}
\end{lem}

\begin{proof}
(a) Since $w_n w_{n+1} \cdots w_{n+k-1} = 0^k$, we obtain $d_{n+k}(w) \le \binom{n+k}{k} d_n(w)$ by  Lemma \ref{lem:binomialinequality}(a). Thus,
\begin{align*}
r_{n+k} &= \frac{d_{n+k}(w)}{(n+k)!} M^{n+k} \\
&\le \frac{1}{(n+k)!} \binom{n+k}{k} d_n(w) M^{n+k} \\
&= \frac{1}{k!} \cdot \frac{d_n(w)}{n!} M^{n+k} = \frac{M^k}{k!} r_n,
\end{align*}
which is less than or equal to $r_n$ by our assumption that $M^k/k! \le 1$.

(b) We first justify the \textit{``$k$ is large enough''} clause: Since $L < 2/\pi$, we have $\frac{1}{2} \cdot \frac{E_k}{k!} \cdot M^k \sim (2/\pi)^{k+1} / L^k \to \infty$ as $k \to \infty$. Thus the inequality $\frac{1}{2} \cdot \frac{E_k}{k!} \cdot M^k \ge 1$ holds for all sufficiently large $k$.

Since $w_n w_{n+1} \cdots w_{n+k-1} = (10)^{k/2}$, we obtain $d_{n+k}(w) \ge \frac{1}{2} \binom{n+k}{k} E_k d_n(w)$ by Lemma \ref{lem:binomialinequality}(b). Thus,
\begin{align*}
r_{n+k} &= \frac{d_{n+k}(w)}{(n+k)!} M^{n+k} \\
&\ge \frac{1}{2} \cdot \frac{1}{(n+k)!} \binom{n+k}{k} E_k d_n(w) M^{n+k} \\
&= \frac{1}{2} \cdot \frac{E_k}{k!} \cdot \frac{d_n(w)}{n!} \cdot M^{n+k}
= \frac{1}{2} \cdot \frac{E_k}{k!} \cdot M^k \cdot r_n,
\end{align*}
which is greater than or equal to $r_n$ by our assumption that $\frac{1}{2} \cdot \frac{E_k}{k!} \cdot M^k \ge 1$.
\end{proof}

% \begin{lem} \label{lem:nalpha}
% If $n_1, n_2, \ldots, n_\alpha, \ldots$ are the values of $n$ at which the state $q$ is changed in the algorithm defined above, then the sequence $\{n_\alpha\}$ is infinite.
% \end{lem}

% \begin{proof}
% Otherwise, there is a last $n$ at which $q$ changes, say $n = N$. If $q$ changes to $0$ at $N$, then $w$ ends with $0^\omega$, so $\gr d_n(w) = 0$. But $L > 0$, so $r_n \to 0$, and so there exists $n > N$ such that $r_n < 1$. But this would change $q$ from $0$ to $10$ in step (1) of the algorithm, a contradiction.

% Similarly, if $q$ changes to $10$ at $N$, then $w$ ends with $(10)^\omega$, so $\gr d_n(w) = 2/\pi$. But $L < 2/\pi$, so $r_n \to \infty$, and so there exists $n > N$ such that $r_n > 1$. But this would change $q$ from $10$ to $0$ in step (1) of the algorithm, a contradiction.
% \end{proof}
We now use the previous lemmas to prove Theorem~\ref{thm:properexists}.
\begin{proof}[Proof of Theorem \ref{thm:properexists}]
In light of Lemma \ref{lem:0^k-inequality}, fix a constant $K$ such that for all $k \ge K$ we have $M^k/k! \le 1$ and $\frac12 \cdot \frac{E_k}{k!} \cdot M^k \ge 1$. Let $n_1, n_2, \ldots, n_\alpha, \ldots$ be the values of $n$ at which the state $q$ is changed in the algorithm defined above, either from $10$ to $0$ or from $0$ to $10$. (It can be proved that the sequence $\{n_\alpha\}$ is infinite, although we do not need this fact in our proof.)
% By Lemma \ref{lem:nalpha}, the sequence $\{n_\alpha\}$ is infinite.

By Remark \ref{rmk:suffices-to-show}, it suffices to show that $\lim_{n\to\infty} {r_n}^{1/n} = 1$. We will do this by showing that
\[ f(n) \le {r_n}^{1/n} \le g(n), \]
where
\[ f(n) = 1/{\left[(n+K)^K \cdot n \right]}^{1/n}
\quad \text{and} \quad
g(n) = {\left(M^{K+2}\right)}^{1/n}. \]
Let $n\ge 1$. We consider four cases:

\begin{itemize}
\item Case (i): $n = n_\alpha$, and $q$ changes from $10$ to $0$ at $n_\alpha$. This means that $r_{n_\alpha} > 1$ but $r_{n_\alpha - 2} \le 1$. Then $r_{n_\alpha} \le M^2 r_{n_\alpha - 2} \le M^2$ (the first inequality uses Lemma \ref{lem:inequality-induction} and the second inequality uses $r_{n_\alpha - 2} \le 1$). Thus
\begin{equation} \label{eqn:nalpha-upperbound}
r_{n_\alpha} \le M^2,
\end{equation}
and so
\[ f(n_\alpha) \le 1 \le {r_{n_\alpha}}^{1/{n_\alpha}} \le {(M^2)}^{1/n_\alpha} \le g(n_\alpha). \]

\item Case (ii): $n_\alpha < n < n_{\alpha+1}$, and $q$ changed from $10$ to $0$ at $n_\alpha$. Note that $q$ remains $0$ from $n_\alpha$ to $n$, so we have $r_{n_\alpha} \ge 1$ and $r_n \ge 1$ and $w_{n_\alpha} \cdots w_n = 0^k$ where $k = n-n_\alpha+1$. If $k \ge K$, then Lemma \ref{lem:0^k-inequality}(a) applies. Thus we need two subcases:
\begin{enumerate}
\item If $k < K$, then
$r_n = r_{n_\alpha+k}
\le M^k r_{n_\alpha}
\le M^k M^2
\le M^{K+2}$
(the first inequality uses Lemma \ref{lem:inequality-induction} and the second inequality uses \eqref{eqn:nalpha-upperbound}).
\item If $k \ge K$, then
$r_n \le r_{n_\alpha}
\le M^2
\le M^{K+2}$
(the first inequality uses Lemma \ref{lem:0^k-inequality}(a) and the second inequality uses \eqref{eqn:nalpha-upperbound}).
\end{enumerate}
In both subcases, we obtain
\[ f(n) \le 1 \le {r_n}^{1/n} \le {(M^{K+2})}^{1/n} = g(n). \]

\item Case (iii): $n = n_\alpha$, and $q$ changes from $0$ to $10$ at $n_\alpha$. This means that $r_{n_\alpha} < 1$ but $r_{n_\alpha - 1} \ge 1$. Then $r_{n_\alpha} \ge \frac{M}{n_\alpha} r_{n_\alpha - 1} \ge \frac{M}{n_\alpha}$ (the first inequality uses Lemma \ref{lem:r_n-inequality} and the second inequality uses $r_{n_\alpha - 1} \ge 1$). Thus
\begin{equation} \label{eqn:nalpha-lowerbound}
r_{n_\alpha} \ge \frac{M}{n_\alpha},
\end{equation}
and so
\[ f(n_\alpha) \le \frac{M}{n_\alpha} \le {r_{n_\alpha}}^{1/n_\alpha} \le 1 \le g(n_\alpha). \]

\item Case (iv): $n_\alpha < n < n_{\alpha+1}$, and $q$ changed from $0$ to $10$ at $n_\alpha$. Note that $q$ remains $10$ from $n_\alpha$ to $n$, so we have $r_{n_\alpha} \le 1$ and $r_n \le 1$ and $w_{n_\alpha} \cdots w_n = (10)^k$ where $k = n-n_\alpha+1$. If $k \ge K$, then Lemma \ref{lem:0^k-inequality}(b) applies. Thus we again need two subcases:
\begin{enumerate}
\item If $k < K$, then
\begin{align*}
r_n = r_{n_\alpha+k}
&\ge \frac{M^k}{(n_\alpha+k)^k} r_{n_\alpha}
\ge \frac{M^k}{(n_\alpha+k)^k} \cdot \frac{M}{n_\alpha} \\
&\ge \frac{1}{(n_\alpha+k)^k} \cdot \frac{1}{n_\alpha}
\ge \frac{1}{(n+K)^K \cdot n}
\end{align*}
(the first inequality uses Lemma \ref{lem:inequality-induction} and the second inequality uses \eqref{eqn:nalpha-lowerbound}).
\item If $k \ge K$, then
\[
r_n \ge r_{n_\alpha}
\ge \frac{M}{n_\alpha}
\ge \frac{1}{n_\alpha}
\ge \frac{1}{(n+K)^K \cdot n}
\]
(the first inequality uses Lemma \ref{lem:0^k-inequality}(b) and the second inequality uses \eqref{eqn:nalpha-lowerbound}).
\end{enumerate}
In both subcases, we obtain
\[ f(n) = \frac{1}{(n+K)^K \cdot n} \le {r_n}^{1/n} \le 1 \le g(n). \]
\end{itemize}

Therefore, in all four cases, $f(n) \le {r_n}^{1/n} \le g(n)$. Therefore, $\lim_{n\to\infty} {r_n}^{1/n} = 1$, because $\lim_{n\to\infty} f(n) = \lim_{n\to\infty} g(n) = 1$.
\end{proof}

By modifying Algorithm 1, we obtain a more general version of Theorem \ref{thm:properexists} which allows different lower and upper growth rates. Let $L, L' \in [0, 2/\pi]$ such that $L \le L'$. We modify the algorithm as follows. Define $M' = (L')^{-1}$, and define $r'_n = \frac{d_n(w)}{n!} (M')^n$. When $q = 10$, instead of checking the condition $r_n > 1$, the algorithm should check the condition $r'_n > 1$ (and set $q = 0$ when this condition is met). The same reasoning as in the proof of Theorem \ref{thm:properexists} now shows that $\liminf_{n\to\infty} r_n^{1/n} = 1$ and $\limsup_{n\to\infty} {(r'_n)}^{1/n} = 1$, from which we conclude the following theorem:
\begin{thm}\label{thm:infsup}
For every $L, L' \in [0, 2/\pi]$ such that $L \le L'$, there exists $w \in \infwords$ such that $\logr d_n(w) = L$ and $\upgr d_n(w) = L'$.
\end{thm}

\section{Peaks}\label{sec:peaks}

This section is dedicated to proving Theorem~\ref{thm:peakgrowthrate}. As we will see, the key constructions in this section are afforded by main theorems in \cite{BFT}. We start by recalling phenomena about finite peak sets for permutations in $\Sn$. Throughout, we reserve the variable $S$ rather than the variable $I$ for \emph{finite} subsets of $\mathbb{N}$. In that light,
given a positive integer $n$ and $S \subseteq [n]$, let $\P$ be the subset of permutations in $\Sn$ with peak set precisely $S$ (as introduced in \cite{BBS}). In other words, $\P=\{\pi \in \Sn \colon \pi_{i-1}<\pi_i>\pi_{i+1} \mbox{ if and only if } i \in S\}$. When $\P \neq 0$ we say $S$ is \emph{admissible}, otherwise $S$ is said to be inadmissible. So for instance, if $S$ has two consecutive positive integers in it, or if $1 \in S$ or $n \in S$, then $S$ is \emph{inadmissible}. 

Our main strategy for computing $|\P|$ is to strategically and iteratively split $S$ into disjoint sets $S=S_L \cup S_R$ with
\[ S_L =
\{i_1,i_2,\ldots,i_{\ell}\}, \ S_R=\{m+1+j_1,m+1+j_2,\ldots,m+1+j_r\}
\]
where $i_1 < i_2 < \cdots < i_{\ell} = m$ and $2=j_1<j_2<\cdots<j_r$. In this case we can compute $|\P|$ in terms of $|\mathcal{P}_{S_L}(m+1)|$ and $|\mathcal{P}_{S_R}(n-(m+1))|$.

\begin{thm}\label{thm:split}\cite[Thm.\@ 3.15]{BFT}
Let $S \subseteq \mathbb{N}$ be admissible, and suppose $S$ can be decomposed as 
\[
S=\{i_1<i_2<\ldots<i_{\ell}=m\} \cup \{m+1+j_1,m+1+j_2,\ldots,m+1+j_r\}.
\]
If $S_L=\{i_1<i_2<\cdots<i_{\ell}=m\}$ and $S_R=\{j_1=2<j_2<\cdots<j_r\}$,
then 
\[
|\P| = \binom{n}{m+1} |\mathcal{P}_{S_L}(m+1)| \cdot |\mathcal{P}_{S_R}(n-(m+1))|.
\]
\end{thm}
Theorem~\ref{thm:split} applies when $S$ has what is referred to as a ``gap of three'', the name coming from the fact that the largest element in $S_L$ and the smallest in $S_R+(m+1)=\{m+1+j_1,m+1+j_2,\ldots,m+1+j_r\}$ differ by exactly three.

In our applications, we use Theorem~\ref{thm:split} to split $S$ into sets of two types, of which one type consists of what \cite{BFT} called \emph{maximal alternating subsets}: subsets of $S$ of the form $\{\ell,\ell+2,\ell+4,\ldots,\ell+2k\}$ with neither $\ell-2$ nor $\ell+2k+2$ in $S$. The naming comes from the fact that a permutation with such a peak set $S$ is alternating on the indices in the interval $[\ell,\ell+2k]$. The relevance of such sets motivates the following: let $\A(S)$ be the partition of $S$ into maximal alternating subsets. For instance if 
\[
S=\{3,6,9,12,14,16,19,21,24,27,29,31,34,36\}
\]
then 
\[
\A(S)=\{\{3\},\{6\},\{9\},\{12,14,16\},\{19,21\},\{24\},\{27,29,31\},\{34,36\}\}.
\]
As introduced in \cite{BFT}, for an admissible set $S \subseteq [n]$ define the statistic
\[\Q_S(n):= \bigcup_{T \supseteq S} \mathcal{P}_T(n) ; \]
that is, the set of permutations in $\Sn$ whose peak set contains $S$.  The following lemma recovers $|\Q_S(n)|$ from $\A(S)$ and represents it in terms of Euler numbers:
 
\begin{lem}\cite[Lem.\@ 5.1]{BFT}\label{lem:Q}\
Let $S \subseteq \mathbb{N}$ be admissible and $m=\max(S)$. For $n \geq m+1$,
\[
|\Q_S(n)| = n! \prod_{A \in \A(S)} \frac{E_{2|A|+1}}{(2|A|+1)!}.
\]
\end{lem}

This can be used together with inclusion/exclusion to establish a formula for $|\P|$ using the formula for $|\Q_S(n)|$:

\begin{lem}\label{lem:inclusion}\cite[Lem.\@ 5.2]{BFT}
Let $S \subseteq \mathbb{N}$ be admissible and $m=\max(S)$. For $n \geq m+1$,
\[
|\P| = \sum_{T \supseteq S} (-1)^{|T-S|}|\Q_T(n)|.
\]
\end{lem}
Notice as a consequence that if $S$ is not a proper subset of an admissible set, then $|\P|=|\Q_S(n)|$. 

To prove Theorem~\ref{thm:peakgrowthrate}, we construct carefully chosen periodic words of the form  $(01(001)^a0^b)^{\omega}$ for $a,b \geq 2$ and show they have growth rates that form a dense set in $[0,1/\sqrt[3]{3}]$. Roughly speaking, the larger $b$ is compared to $a$, the closer to $0$ the growth rate is, whereas if $a$ is large compared to $b$, the growth rates gets closer to $1/\sqrt[3]{3}$. Some care is needed to balance the ratio $a/b$ to allows us to get arbitrarily close to growth rates in the interval $[0,1/\sqrt[3]{3}]$. This is done by first finding an explicit formula for the growth rates of the aforementioned words. To do so, we first identify $w|_{n-1}$ with the subset $S \subseteq [n-1]$ that contains positions where $w|_{n-1}$ is $1$. Then, we partition $[n-1]$ using Theorem~\ref{thm:split} into two types of subsets: (1) nearly maximal alternating subsets as discussed above and (2) subsets that are large but with only two elements of $S$. Theorem~\ref{thm:split} then applies iteratively because these sets $S$ coming from the words $(01(001)^a0^b)^{\omega}$ have many gaps of three. To count the permutations associated to sets of type (1) we use Lemma~\ref{lem:inclusion}. To count permutations associated to sets of type (2) we use the following:

\begin{lem}\label{lem:piece}
For a positive integer $n \geq 4$ we have $\P=2^{n-3} \cdot (n-4)(n-1)$ if $S=\{2,n-1\}$ and $\P=2^{n-2} \cdot (n-2)$ if $S=\{2\}$.
\end{lem}
\begin{proof}
These follow from Theorems 6 and 7 in \cite{BBS}.
\end{proof}

We can now determine the growth rates of our periodic words:

\begin{thm}\label{thm:exactpeakgrowthrate}
Let $a,b \geq 2$ and $w=(01(001)^a0^b)^{\omega}$. Then  
\[
\gr p_n(w) = \left(\frac{1}{(b+5)!}\right)^{1/((3a-3)+(b+5))} \cdot \left(\frac{1}{3}\right)^{\frac{a-1}{3(a-1)+(b+5)}} \cdot 2^{\frac{b+2}{3(a-1)+b+5}}
\]
\end{thm}
\begin{proof}
First observe that $w$ can be rewritten as 
\[
w=w^{(1)}(w^{(2)}w^{(3)})^{\omega}
\]
where $w^{(1)}=(010)^a, w^{(2)}=(010)0^{b-1}(010),$ and $w^{(3)}=(010)^{a-1}$. The sum of the lengths of the words $w^{(2)}$ and $w^{(3)}$ is $3a+b+2$. For sufficiently large $n$, select (by the Division algorithm) $1 \leq r \leq 3a+b+2$ so that $n-(3a)=q(3a+b+2)+r$ for some positive integer $q$.

First consider if $1 \leq r < b+5$. Then
\[ w|_n = w^{(1)}(w^{(2)}w^{(3)})^q [w^{(2)}|_r]. 
\]
Now notice we can apply Theorem~\ref{thm:split} repeatedly to $p_n(w)$ because, in each instance of $w^{(i)} w^{(j)}$, the last $1$ in $w^{(i)}$ is exactly three positions before the first $1$ in $w^{(j)}$. From this,
\[
p_n(w) = \frac{n!}{(3a)! \cdot (b+5)!^q \cdot (3a-3)!^q \cdot r!} \cdot \mathcal{P}_S(3a) \cdot \mathcal{P}_{S'}(b+5)^q \cdot \mathcal{P}_{S''}(3a-3)^q \cdot \mathcal{P}_{S'''}(r)
\]
where the indexed sets are given by 
\begin{align*}
&S=\{2,5,8,\ldots\} \cap [3a] & &S'=\{2,b+4\} \\
&S''=\{2,5,8,\ldots\} \cap [3a-3] & &S'''=\{2,b+4\} \cap [r-1] = \{2\}
\end{align*}
Lemma~\ref{lem:piece} gives us that $\mathcal{P}_{S'}(b+5)=2^{b+2}(b+1)(b+4)$ and similarly $\mathcal{P}_{S'''}(r)=2^{r-2}(r-2)$. Since $S,S'',S'''$ have no admissible supersets, Lemma~\ref{lem:inclusion} and Lemma~\ref{lem:Q} together give us
\[
\mathcal{P}_S(3a)= (3a)! \left( \frac{E_3}{3!} \right)^a, \ \ \mathcal{P}_{S''}(3a-3)= (3a-3)! \left( \frac{E_3}{3!} \right)^{a-1}, \ \ \mathcal{P}_{S'''}(r) = 2^{r-2}(r-2).
\]
Altogether then, we have $p_n(w)$ is
\scriptsize
\begin{align*}
&\frac{n!}{(3a)! \cdot (b+5)!^q \cdot (3a-3)!^q \cdot r!} 
 \cdot (3a)! \left( \frac{E_3}{3!} \right)^a \cdot (2^{b+2}(b+1)(b+4))^q \cdot \left( (3a-3)! \left( \frac{E_3}{3!} \right)^{a-1} \right)^q \cdot  2^{r-2}(r-2) \\
&= \frac{n!}{(3a)! \cdot (b+5)!^q \cdot (3a-3)!^q \cdot r!}  \cdot (3a)! \cdot (3a-3)!^q \cdot \left(\frac{1}{3}\right)^{a+q(a-1)} \cdot (2^{b+2}(b+1)(b+4))^q \cdot 2^{r-2}(r-2).
\end{align*}
\normalsize
For large $n$ the only terms that contribute to the asymptotics are those with a factor of $q$, and all polynomial factors can be eliminated, so we can approximate $\frac{p_n(w)}{n!}$ arbitrarily close by
\[
\frac{1}{(b+5)!^q} \cdot \frac{1}{3}^{q(a-1)} \cdot 2^{q(b+2)}
\]
and so $\gr p_n(w)$ is approximated arbitrarily close by 
\[
\left(\frac{1}{(b+5)!}\right)^{1/((3a-3)+(b+5))} \cdot \left(\frac{1}{3}\right)^{\frac{a-1}{3(a-1)+(b+5)}} \cdot 2^{\frac{b+2}{3(a-1)+b+5}}
\]
because $n=3a+q((3a-3)+(b+5))+r$ and the $3a$ and $r$ contributions are negligible as $n$ is large.

If instead we have $b+5<r\leq 3a+b+2$, then 
\[
w|_n = w^{(1)}(w^{(2)}w^{(3)})^qw^{(2)} [w^{(3)}|_{r-(b+5)}] 
\]
and so the same analysis gives
\scriptsize
\begin{align*}
p_n(w) &= \frac{n!}{(3a)! \cdot (b+5)!^{q+1} \cdot (3a-3)!^q \cdot (r-(b+5))!} \\
 &\cdot (3a)! \left( \frac{E_3}{3!} \right)^a \cdot (2^{b+2}(b+1)(b+4))^{q+1} \cdot \left( (3a-3)! \left( \frac{E_3}{3!} \right)^{a-1} \right)^q  \cdot (r-(b+5))! \left( \frac{E_3}{3!} \right)^{\left \lfloor \frac{r-(b+5)}{3} \right \rfloor} \\
&= \frac{n!}{(3a)! \cdot (b+5)!^{q+1} \cdot (3a-3)!^q \cdot (r-(b+5))!} \\
 &\cdot (3a)! \cdot (3a-3)!^q \cdot (r-(b+5))! \cdot \left( \frac{1}{3} \right)^{a+q(a-1)+\frac{r-(b+5)}{3}} \cdot (2^{b+2}(b+1)(b+4))^{q+1}
\end{align*}
\normalsize
A similar analysis gives that 
$\gr p_n(w)$ is approximated arbitrarily close by 
\[
\left(\frac{1}{(b+5)!}\right)^{1/((3a-3)+(b+5))} \cdot \left(\frac{1}{3}\right)^{\frac{a-1}{3(a-1)+(b+5)}} \cdot 2^{\frac{b+2}{3(a-1)+b+5}}
\]
\end{proof}
We now prove Theorem~\ref{thm:peakgrowthrate}, showing that for each $L \in [0,1/\sqrt[3]{3}]$ and each $\epsilon > 0$  there is a periodic word with growth rate within $\epsilon$ of $L$. The proof hinges heavily on Theorem~\ref{thm:exactpeakgrowthrate}.

\begin{proof}[Proof of Theorem~\ref{thm:peakgrowthrate}]
Let $L \in [0,1/\sqrt[3]{3}]$ and $\epsilon>0$. For simplicity let $c=1/\sqrt[3]{3}$ and $\gamma = \log(c/L)$. By choice of $L$, we have $\gamma>0$. For a positive integer $m$ let $w$ be the word $(01(001)^a0^b)^{\omega}$ where $a=\lfloor m \log(m) \rfloor$ and $b=\lfloor \gamma m \rfloor$. Note that $w$ depends on $m$; we claim that we can choose $m$ sufficiently large so that $|\gr p_n(w)-L|<\epsilon$. To that end, consider $\log \gr p_n(w)$ which by Theorem~\ref{thm:exactpeakgrowthrate} is given by
\scriptsize
\begin{align*}
\log \gr p_n(w) &= \frac{1}{3\lfloor m \log(m) \rfloor + (\lfloor \gamma m \rfloor + 5)} \\ & \hspace{0.5in} \cdot \left( 3(\lfloor m \log(m) \rfloor -1)\log(c) + (\lfloor \gamma m \rfloor + 2)\log(2) - \log((\lfloor \gamma m \rfloor + 5)!)\right)
\end{align*}
\normalsize
Distributing through the parentheses we get three terms. As $m \to \infty$, the first term approaches $\log(c)$, the second term vanishes, and by Stirling's formula the third term approaches $-\gamma$. Altogether then, as $m \to \infty$, we have $\log \gr p_n(w) \to -\gamma \log(c)$ and so $\gr p_n(w) \to e^{-\gamma \log(c)} = L$, and the result follows.
\end{proof}

\section{Conclusions}\label{sec:conclusion}
We end the article with pertinent questions of interest. The strategy employed in Section~\ref{sec:descents} for constructing a set $I \subseteq \mathbb{N}$ that achieves a particular limit $L \in \left[0,2/\pi \right]$ builds $I$ in chunks of sets of consecutive even integers that are themselves spread out from each other. One might think the same strategy can be employed for proving that limits of growth rates of peak sets are themselves dense in $\left[0,1/\sqrt[3]{3}\right]$, however the issue that arises is that not having gaps of three (see the discussion after Theorem~\ref{thm:split}) makes it difficult to use theorems such as Theorem~\ref{thm:split}. One could try to use Lemma~\ref{lem:inclusion} instead, however for a given $n$, $S=I \cap [n-1]$ might have many admissible supersets, and the inclusion--exclusion gets out of hand quite quickly. This calls for a new method to address the following question:

\begin{prob}
For which $L \in \left[0,1/\sqrt[3]{3}\right]$ is there a set $I \subseteq \mathbb{N}$ such that $\gr p_n(I)=L$?
\end{prob}

Recently, attention has been given to a more general concept of a descent called an $X$-descent. For a set $X \subseteq \mathbb{N} \times \mathbb{N}$, and $I \subseteq \mathbb{N}$, the permutations in $\Sn$ with $X$-descent $I$ is the set
\[
D_{X,n}(I) := \{\pi \in \Sn \colon (\pi_i,\pi_{i+1}) \in X \mbox{ if and only if } i \in I\}.
\]
For instance, if $X=\{(i,j) \colon i>j\}$ then $D_{X,n}(I)$ is the precisely the set of permutations in $\Sn$ with descent set $I$ (see for instance \cite{Grinberg},\cite{Omar} and references therein). Now letting $d_{X,n}(I)=|D_{X,n}(I)|$, what can we say about the growth rates of the sequences $\{d_{X,n}(I)\}$ for various $X$?

\begin{prob}
Determine all $L \in \mathbb{R}$ for which there are sets $X \subseteq \mathbb{N} \times \mathbb{N}$ and $I \subseteq \mathbb{N}$ such that 
\[
\lim_{n \to \infty} \left( d_{X,n}(I)/n! \right)^{1/n}=L.
\]
\end{prob}

Our developments made key insights on the range of growth rates for descents and peaks. In the case of descents, we noted that $\mathbb{E}[d_n(w)]=2\left(1/2\right)^n \cdot n!$ and so a quick calculation shows the expected growth rate is $1/2$; a similar calculation by linearity of expectation establishes an expected value of $1/3$ for peaks. One insight we don't have however is what the distribution of growth rates of descents (and similarly peaks) looks like.

\begin{prob} \label{prob:distribution}
Determine the distribution of the numbers $\gr d_n(I)$ over all $I \subseteq \mathbb{N}$. Similarly, determine the distribution of the numbers $\gr p_n(I)$ over all admissible $I \subseteq \mathbb{N}$.
\end{prob}

For $d_n(I)$, this question has been studied by Luck \cite{Luck} from the perspective of statistical mechanics. In this setting, the permutation is replaced with a sequence of identically distributed random numbers. The probability that the random sequence has an up-down sequence that matches a given word $w$ is denoted by $P_n(w)$, and the central quantity of interest is the decay rate $\alpha$ such that $P_n(w)$ is asymptotically the product of $e^{-\alpha n}$ and some sub-exponential factor. This $\alpha$ can be expressed in terms of the probability $P_n(w)$ as $\alpha = \lim_{n\to\infty} -\ln(P_n(w))/n$, and it can be expressed in terms of our growth rate $\gr d_n(w)$ as $\alpha = -\ln( \gr d_n(w))$. Luck views $\alpha$ as the \emph{embedding entropy} of $w$. In addition to his deep dives into a handful of illuminating examples (see our Remark \ref{rmk:chirping}), Luck also studies the statistics of $P_n(w)$ when $w$ is randomly chosen from various statistical ensembles. For the uniform ensemble, i.e.\@ all $w|_{n-1}$ are equally likely, he finds numerically (non-rigorously) that the ``typical value'' of $\alpha$ is $\alpha_0 \approx 0.806361$, in the sense that the distribution of $\alpha$ is concentrated around $\alpha_0$. Interestingly, $\alpha$ does not have expected value $\ln(2)$ even though $d_n(w)$ has expected value $1/2^{n-1}$; this is an instance of the fact that $\mathbb{E}[\ln(X)]$ and $\ln(\mathbb{E}[X])$ are not necessarily equal. We view our work as running in parallel with Luck's work, taking the combinatorial road while Luck takes the statistical road.

\section*{Acknowledgments}
We thank the two anonymous referees, whose recommendations and suggestions helped us improve our paper. The first author is partially supported by research funds from York University, and NSERC Discovery Grant \#RGPIN-2025-06304. The second author is partially supported by start-up research funds from California State University, Los Angeles.


\begin{thebibliography}{00}

\bibitem{Bencs}
Ferenc Bencs.
\newblock Some coefficient sequences related to the descent polynomial.
\newblock {\em European J. Combin.}, 98:103396, 2021.

\bibitem{BHR}
Edward~A. Bender, William~J. Helton, and L.~Bruce Richmond.
\newblock Asymptotics of permutations with nearly periodic patterns of rises
  and falls.
\newblock {\em Electron. J. Combin.}, 10:{\#}R40, 2003.

\bibitem{BBS}
Sara Billey, Krzysztof Burdzy, and Bruce~E. Sagan.
\newblock Permutations with given peak set.
\newblock {\em J. Integer Seq.}, 16:{ }Article 13.6.1, 2013.

\bibitem{BFT}
Sara Billey, Matthew Fahrbach, and Alan Talmage.
\newblock Coefficients and roots of peak polynomials.
\newblock {\em Experiment. Math.}, 25(2):165--175, 2016.

\bibitem{BFW}
Francis C.~S. Brown, Thomas M.~A. Fink, and Karen Willbrand.
\newblock On arithmetic and asymptotic properties of up--down numbers.
\newblock {\em Discrete Math.}, 307(14):1722--1736, 2007.

\bibitem{deBruijn}
Nicolaas Govert~de Bruijn.
\newblock Permutations with given ups and downs.
\newblock {\em Nieuw Arch. Wiskd.}, 18:61--65, 1970.

\bibitem{DHIO}
Alexander Diaz-Lopez, Pamela~E. Harris, Erik Insko, and Mohamed Omar.
\newblock A proof of the peak polynomial positivity conjecture.
\newblock {\em J. Combin. Theory Ser. A}, 149:21--29, 2017.

\bibitem{DHIOS}
Alexander Diaz-Lopez, Pamela~E. Harris, Erik Insko, Mohamed Omar, and Bruce~E.
  Sagan.
\newblock Descent polynomials.
\newblock {\em Discrete Math.}, 342(6):1674--1686, 2019.

\bibitem{ET}
Sergi Elizalde and Justin~M. Troyka.
\newblock Exact and asymptotic enumeration of cyclic permutations according to
  descent set.
\newblock {\em J. Combin. Theory Ser. A}, 165:360--391, 2019.

\bibitem{Grinberg}
Darij Grinberg and Richard~P. Stanley.
\newblock The {Redei}--{Berge} symmetric function of a directed graph.
\newblock {\em \textup{\texttt{arXiv:2307.05569}}}, 2023.

\bibitem{KantarciOguz}
Ezgi Kantarc{\i}~O{\u{g}}uz.
\newblock Connecting descent and peak polynomials.
\newblock {\em Hacet. J. Math. Stat.}, 53(2):488--494, 2024.

\bibitem{Luck}
Jean-Marc Luck.
\newblock On the frequencies of patterns of rises and falls.
\newblock {\em Phys. A}, 407:252--275, 2014.

\bibitem{Macmahon}
Percy~Alexander MacMahon.
\newblock {\em Combinatory Analysis}, volume~1.
\newblock Cambridge University Press, 1915.

\bibitem{Niven}
Ivan Niven.
\newblock A combinatorial problem of finite sequences.
\newblock {\em Nieuw Arch. Wiskd.}, 16:116--123, 1968.

\bibitem{Omar}
Mohamed Omar.
\newblock Permutations with a given {$X$}-descent set.
\newblock {\em \textup{\texttt{arXiv:2402.10443}}}, 2024.

\bibitem{Stanley}
Richard~P. Stanley.
\newblock {\em Enumerative Combinatorics, Volume 1}, volume~49 of {\em
  Cambridge Stud. Adv. Math.}
\newblock Cambridge Univ. Press, 2nd edition, 2012.

\bibitem{vatter2015permutation}
Vincent Vatter.
\newblock Permutation classes.
\newblock {\em Handbook of enumerative combinatorics}, pages 754--833, 2015.

\end{thebibliography}
\end{document}